\documentclass{iopart}
\usepackage{iopams}
\eqnobysec
\usepackage{graphicx}
\usepackage{epstopdf}
\usepackage{mathrsfs}

\newtheorem{lemma}{Lemma}[section]
\newtheorem{theorem}[lemma]{Theorem}
\newtheorem{remark}{Remark}[section]
\newtheorem{Rule}{Rule}[section]

\newtheorem{example}[lemma]{Example}
\newtheorem{corollary}[lemma]{Corollary}

\newtheorem{assumption}{Assumption}[section]
\newcommand{\proof} [1]
   { {\bf Proof.} #1 \hfill\opensquare \\}

\def\p{\partial}
\def\d{\delta}

\def\ep{\varepsilon}
\def\l{\langle}
\def\r{\rangle}

\def\R{\mathcal R}
\def\D{\mathscr D}
\def\X{X}
\def\Y{Y}

\def\a{\alpha}

\def\d{\delta}

\begin{document}
\jl{5}
\title[]
{Hanke-Raus heuristic rule for variational regularization in Banach spaces}

\author{Qinian Jin}

\medskip

\address{Mathematical Sciences Institute, Australian National University,
Canberra, ACT 0200, Australia}


\eads{Qinian.Jin@anu.edu.au}



\begin{abstract}
We generalize the heuristic parameter choice rule of Hanke-Raus for quadratic regularization to general
variational regularization for solving linear as well as nonlinear ill-posed inverse problems in Banach spaces.
Under source conditions formulated as variational inequalities, we obtain a posteriori error estimates in
term of Bregman distance. By imposing certain conditions on the random noise, we establish four convergence
results; one relies on the source conditions and the other three do not depend on any source conditions.
Numerical results are presented to illustrate the performance.
\end{abstract}





%

\section{Introduction}

Inverse problems frequently occur in many practical applications in natural sciences, engineering and medicine
whenever one searches for unknown causes based on observations of their effects (see \cite{EHN1996,SKHK2012}).
In this paper we consider inverse problems of the form
\begin{equation}\label{1}
F(x) =y,
\end{equation}
where $F: \D(F) \subset \X \to \Y$ is an operator between two Banach spaces $\X$ and $\Y$ with domain $\D(F)$.
The norms in $\X$ and $\Y$ are always denoted by the same notation $\|\cdot\|$ which should be clear from the context.

Throughout the paper we assume that (\ref{1}) has a solution. In general (\ref{1}) may have many solutions.
In order to find the one with the desired feature, we choose a proper, lower semi-continuous, convex
function $\R: \X \to [0, \infty]$ and determine a solution $x^\dag$ with the property
\begin{equation*}
\R(x^\dag) = \min\left\{\R(x):  x\in \D(F) \mbox{ and } F(x) = y\right\}
\end{equation*}
which is called an $\R$-minimizing solution of (\ref{1}). Because of measurement errors in practical applications,
instead of $y$ we only have a noisy data $\tilde{y}$ whose noise level is denoted by
\begin{equation*}
\d:=\|\tilde{y}-y\|.
\end{equation*}
Due to the inherent ill-posedness of inverse problems, the computation of $x^\dag$ from $\tilde{y}$ requires a
regularization method. Variational regularization is a family of prominent methods in which a minimizer
\begin{equation}\label{Tik}
\tilde{x}_\a \in \arg\min_{x\in \D(F)} \left\{{\mathcal T}_\a(x):=\|F(x) -\tilde{y}\|^r + \a \R(x)\right\}
\end{equation}
with $1< r<\infty$ is used to approximate $x^\dag$. Here, the regularization parameter $\a>0$ plays a crucial role
for the accuracy of approximation. How to choose $\a$ is indeed the most challenging and most important
question for variational regularization.

Many parameter choice rules have been proposed to choose the regularization parameter in either
{\it a priori} or {\it a posteriori} ways, including the famous discrepancy principle and its variants
(\cite{EHN1996,EKN1989,HM2012,Jin1999,JH1999,SEK1993,TJ2003}).  All these rules require accurate knowledge
of the noise level $\d$ to obtain satisfactory approximate solutions. In real world applications such
noise level information is not always available or reliable. Overestimation or underestimation on noise level
may lead to a significant loss of accuracy when using these rules. It is therefore necessary
to consider purely data driven parameter choice rules that avoid using knowledge of noise level.

For quadratic regularization in Hilbert spaces several heuristic parameter choice rules not using information
on the noise level have been proposed, including the generalized cross validation \cite{W1977},
the $L$-curve method \cite{HL1993}, the quasi-optimality criterion \cite{KN2008,TGK1979} and
the Hanke-Raus rule \cite{EHN1996,HR1996}. In this paper we will extend the Hanke-Raus rule to study
variational regularization for linear as well as nonlinear inverse problems in Banach spaces.

When the approximation error between $\tilde{x}_\a$ and $x^\dag$ is measured by a certain ``metric" $D(\cdot, \cdot)$,
a satisfactory choice of $\a$ should make $D(\tilde{x}_\a, x^\dag)$ as small as possible. This function however is 
not computable because it involves the sought solution $x^\dag$. The basic idea behind the Hanke-Raus rule is to 
find a computable surrogate $\Theta(\a, \tilde{y})$ such that $\Theta(\a, \tilde{y})$ and $D(\tilde{x}_\a, x^\dag)$ 
have the similar sharp upper bounds in the worst case scenario when the sought
solution satisfies the standard source conditions and to choose the regularization parameter to be a global minimizer of
$\a\to \Theta(\a, \tilde{y})$ over a certain interval $(0, \a_0]$, where $\a_0>0$ is a given number. One may refer to 
\cite{EHN1996,HR1996} for detailed explanations. For the conventional
quadratic regularization in Hilbert spaces which corresponds to (\ref{Tik}) with $r=2$, $F$ a bounded
linear operator and $\R(x) = \|x\|^2$, several choices of $\Theta$ were proposed in \cite{EHN1996,HR1996};
in particular, the function
\begin{equation*}
\Theta(\a, \tilde{y}) := \frac{\|F(\tilde{x}_\a)-\tilde{y}\|^2}{\a}
\end{equation*}
was considered, see \cite[\S 4.5]{EHN1996}. Although Bakushinskii's veto \cite{B1984} says that a heuristic
parameter choice rule does not lead to a convergent regularization method for ill-posed inverse problems
in the worst case scenario, partial theoretical justification of the use of this parameter choice rule was provided
in \cite{EHN1996,HR1996} where a posteriori error estimates were derived under the source conditions
$x^\dag \in {\mathscr R}(F^*F)^\nu)$ with $\nu>0$ and a convergence result was established under additional conditions
on the randomness of noise.  The parameter choice rule of Hanke-Raus was recently extended to study the convex
variational regularization (\cite{JL2010}) which corresponds to (\ref{Tik}) with $r=2$ and $\Y$ a Hilbert space
and the constrained nonlinear Tikhonov regularization in Hilbert spaces (\cite{IJ2011}) which corresponds to (\ref{Tik})
with $r=2$, $\X$ and $\Y$ being Hilbert spaces and $\R(x) = \|x\|^2 +\iota_{\mathcal C}(x)$,
where ${\mathcal C}\subset \X$ is a closed convex set representing constraints on solutions and $\iota_{\mathcal C}$
denotes the indicator function of ${\mathcal C}$. The theoretical results in \cite{IJ2011,JL2010} were obtained
under the source condition
\begin{equation}\label{sc}
{\mathscr R}(F'(x^\dag)^*) \cap \p \R(x^\dag) \ne \emptyset
\end{equation}
on the unknown solution $x^\dag$, where $F'(x^\dag)$ denotes the Fr\'{e}chet derivative of $F$ at $x^\dag$
in case $F$ is Fr\'{e}chet differentiable and $\p \R$ denotes the subdifferential of $\R$. This source condition
is restrictive and is difficult to check in practical applications,

The following questions arise naturally: Is it possible to generalize Hanke-Raus rule to study the general
variational regularization (\ref{Tik}) in Banach spaces? If yes, can we derive the correpsonding a posteriori error
estimates under general source conditions? Can we prove convergence results under certain conditions on the randomness
of noise without using any source conditions on the unknown solution? In this paper we will give affirmative answers
to the above questions. A natural formulation of Hanke-Raus rule
in the context of (\ref{Tik}) is to choose the regularization parameter $\a^*\in (0, \a_0]$ such that
\begin{equation*}
\Theta(\a_*, \tilde{y})= \min_{\a\in (0, \a_0]} \left\{\Theta(\a,\tilde{y}) := \frac{\|F(\tilde{x}_\a)-\tilde{y}\|^r}{\a}\right\}.
\end{equation*}
It should be point out that finding a global minimizer of this $\Theta(\a, \tilde{y})$ over $(0, \a_0]$ is highly
nontrivial and could be very time-consuming. For the purpose of numerical implementation, we restrict the search
of a global minimizer of $\Theta(\a, \tilde{y})$ to a discrete exponential grid. This leads us to propose the following
version of Hanke-Raus parameter choice rule.

\begin{Rule}\label{HRR}
Let $\a_0>0$ and $0<q<1$ be given numbers and set
\begin{equation*}
\Delta_q = \{\a_0 q^j: j=0, 1, \cdots\}.
\end{equation*}
We then define $\a_*: =\a_*(\tilde{y}) \in \Delta_q$ such that
\begin{equation*}
\a_* \in \arg \min_{\a\in \Delta_q} \left\{ \Theta(\a, \tilde{y}):= \frac{\|F(\tilde{x}_\a)-\tilde{y}\|^r}{\a}\right\}.
\end{equation*}
\end{Rule}

The number $\a_0$ in Rule \ref{HRR} is preassigned. One can use the minimizing property of
$\tilde{x}_\a$ to show that $\|F(\tilde{x}_\a)-\tilde{y}\|^r/\a\rightarrow 0$ as $\a \rightarrow \infty$.
Thus, if $\a_0$ is too large, it is very likely that Rule 1.1 will output a large number $\a_*$ and hence
result in an approximate solution with large error. If $\a_0$ is too small, the
resulting approximate solution is too oscillatory to give information on the sought solution.
The choice of $\a_0$ usually depends on a rough guess of the optimal regularization parameter.

In this paper we will provide theoretical justifications on the use of Rule \ref{HRR} in practical applications,
In Section 2 we will derive some a posteriori error estimates under source conditions formulated as
variational inequalities and in Section 3 we will establish various convergence results under certain
conditions on the random noise using or without using source conditions on the sought solutions. We will
provide numerical results in Section 4 to illustrate the performance of Rule \ref{HRR}.

We conclude this section by collecting notation and terminology that will be used. Given a Banach space $\X$ we
use $\X^*$ to denote its dual space. The duality pairing between $\X$ and $\X^*$ is denoted by $\l \cdot, \cdot\r$.
The weak convergence and strong convergence are denoted by $``\rightharpoonup"$ and $``\rightarrow"$ respectively.
For a bounded linear operator $A: \X \to \Y$ between Banach spaces, we use $A^*:\Y^*\to \X^*$ to denote its
adjoint. We also use ${\mathscr N}(A)$ and ${\mathscr R}(A)$ to denote the null space and range space of $A$
respectively. When $\X$ is reflexive, the annihilator of ${\mathscr N}(A)$ equals the closure of ${\mathscr R}(A^*)$ in
$\X^*$, i.e.
\begin{equation*}
{\mathscr N}(A)^\perp:= \{\xi\in \X^*: \l \xi, x\r =0 \, \ \forall x\in {\mathscr N}(A)\} = \overline{{\mathscr R}(A^*)}.
\end{equation*}
For a proper convex function $\R: \X\to [0, \infty]$, we denote by $\p \R$ its subdifferential, i.e.
\begin{equation*}
\p \R(x) =\{ \xi\in \X^*: \R(\bar x) \ge \R(x) + \l \xi, \bar x-x\r \, \  \forall \bar x \in \X\}, \quad x\in \X.
\end{equation*}
Given $\xi \in \p \R(x)$ we define
\begin{equation*}
D_\xi \R(\bar x, x) := \R(\bar x) -\R(x) -\l \xi, \bar x-x\r, \quad \bar x \in \X
\end{equation*}
which is called the Bregman distance induced by $\R$ at $x$ in the direction $\xi$.

Throughout this paper  we always assume that $\X$ and $\Y$ are reflexive, $F$ is weakly closed, $\R$ is proper, lower semi-continuous
and convex, (\ref{1}) has a solution in $\D(\R)$, and ${\mathcal T}_\a$ is coercive for every $\a>0$. These
conditions guarantee that (\ref{1}) has an $\R$-minimizing solution and (\ref{Tik}) has a minimizer
$\tilde{x}_\a$ for every $\a>0$.

\section{\bf A posteriori error estimates}
\setcounter{equation}{0}

In this section we will derive {\it a posteriori} error estimates on $\tilde{x}_{\a_*}$ with $\a_*$ chosen by
Rule \ref{HRR} under the following source conditions on an $\R$-minimizing solution $x^\dag$
of (\ref{1}) formulated as variational inequalities, where ${\mathcal M}_\rho:=\{ x\in \D(F): \R(x)<\rho\}$.

\begin{assumption}\label{A1}
{\it $\p \R(x^\dag) \ne \emptyset$ and there exist $\xi^\dag\in \p \R(x^\dag)$, $0\le \beta<1$ and a concave index
function $\varphi: [0, \infty)\to [0, \infty)$ such that
\begin{equation}\label{3.0}
\l \xi^\dag, x^\dag-x\r \le \beta D_{\xi^\dag}\R(x, x^\dag) + \varphi(\|F(x) -F(x^\dag)\|)
\end{equation}
for all $x\in {\mathcal M}_\rho$ with $\rho >\R(x^\dag)$. Here $\varphi$ is called an
index function if it is continuous and strictly increasing with $\varphi(0) = 0$.
}
\end{assumption}

Assumption \ref{A1} combines the smoothness properties of solutions and the structural conditions of
the nonlinear operator into a single condition, unlike the traditional treatment in which smoothness
conditions and nonlinearity conditions are separated. This source condition with $\varphi(t)= Ct$ was
first introduced in \cite{HKPS2007} for the derivation of convergence rates for nonlinear Tikhonov
regularization in Banach spaces. Its general form was used later, see \cite{HM2012,HY2010} for instance.
One may refer to \cite{HKPS2007,HY2010,SKHK2012} for detailed discussions, including various specific
source conditions that imply Assumption \ref{A1}.

In deriving the a posteriori error estimate under Assumption \ref{A1}, we will use the function
\begin{equation}\label{Phi}
\Phi(t) := \frac{t^r}{\varphi(t)},  \quad t>0.
\end{equation}
Since $\varphi$ is a concave index function and $r>1$, $\Phi$ is also an index function and its inverse
$\Phi^{-1}: (0, \infty) \to (0, \infty)$ is well-defined (\cite{HM2012}).

\begin{theorem}\label{thm1}
Let $x^\dag$ be an $\R$-minimizing solution of (\ref{1}) satisfying Assumption \ref{A1} and let $\a_*\in \Delta_q$
be determined by Rule \ref{HRR}. If $\tilde{x}_{\a_*} \in {\mathcal M}_\rho$ and $\d_*:= \|F(\tilde{x}_{\a_*})-\tilde{y}\|\ne 0$, then there holds
\begin{equation*}
D_{\xi^\dag}\R(\tilde{x}_{\a_*}, x^\dag) \le C\left(1+\frac{\d^r}{\d_*^r}\right) \left(\d^r + \varphi(\d +\d_*)\right),
\end{equation*}
where $\d= \|y-\tilde{y}\|$ is the noise level and $C$ is a constant depending only on $\a_0$, $q$, $r$ and $\beta$.
\end{theorem}

\proof{
We first claim that if $\tilde{x}_\a \in {\mathcal M}_\rho$ then
\begin{equation}\label{3.1}
D_{\xi^\dag} \R(\tilde{x}_\a, x^\dag) \le \frac{1}{1-\beta} \left(\frac{\d^r}{\a} + \varphi\left(\d + \|F(\tilde{x}_\a)-\tilde{y}\|\right)\right)
\end{equation}
and
\begin{equation}\label{3.2}
\|F(\tilde{x}_\a)-\tilde{y}\| \le 5 \d + \Phi^{-1} (2^r \a).
\end{equation}
To see this, by using the minimizing property of $\tilde{x}_\a$ we have
\begin{equation*}
\|F(\tilde{x}_\a)-\tilde{y}\|^r +\a \R(\tilde{x}_\a)
\le \|y - \tilde{y}\|^r + \a \R(x^\dag).
\end{equation*}
In view of the definition of the Bregman distance, this gives
\begin{equation*}
\|F(\tilde{x}_\a)-\tilde{y}\|^r + \a D_{\xi^\dag} \R(\tilde{x}_\a, x^\dag)
\le \d^r + \a \l \xi^\dag, x^\dag-\tilde{x}_\a\r.
\end{equation*}
By virtue of Assumption \ref{A1} we further have
\begin{eqnarray*}
\|F(\tilde{x}_\a)-\tilde{y}\|^r + \a D_{\xi^\dag} \R(\tilde{x}_\a, x^\dag)
& \le \d^r + \a\beta D_{\xi^\dag}\R(\tilde{x}_\a, x^\dag)\\
& \quad \, + \a \varphi(\|F(\tilde{x}_\a)-y\|).
\end{eqnarray*}
Because $0\le \beta<1$, we therefore obtain (\ref{3.1}) and
\begin{equation*}
\|F(\tilde{x}_\a) - \tilde{y}\|^r \le \d^r + \a \varphi\left(\|F(\tilde{x}_\a)-y\|\right).
\end{equation*}
By using the inequality $(a+b)^r \le 2^{r-1} (a^r + b^r)$ for $a, b\ge 0$, we obtain
\begin{equation*}
\|F(\tilde{x}_\a)-y\|^r \le 2^r \d^r + 2^{r-1} \a \varphi(\|F(\tilde{x}_\a)-y\|).
\end{equation*}
If $2^r \d^r \ge 2^{r-1} \a\varphi(\|F(\tilde{x}_\a)-y\|)$, we then obtain
\begin{equation*}
\|F(\tilde{x}_\a)-y\| \le 2^{1+1/r}\d \le 4 \d;
\end{equation*}
if $2^r \d^r < 2^{r-1} \a\varphi(\|F(\tilde{x}_\a)-y\|)$, we have
\begin{equation*}
\|F(\tilde{x}_\a)-y\|^r \le 2^r \a \varphi(\|F(\tilde{x}_\a)-y\|)
\end{equation*}
which shows that $\Phi(\|F(\tilde{x}_\a)-y\|) \le 2^r \a$ and hence $\|F(\tilde{x}_\a)-y\| \le \Phi^{-1}(2^r\a)$.
Combining the estimates from the two cases we thus obtain (\ref{3.2}).

Since we have assumed $\tilde{x}_{\a_*} \in {\mathcal M}_\rho$, we may use (\ref{3.1}) to derive that
\begin{equation}\label{3.3}
D_{\xi^\dag} \R(\tilde{x}_{\a_*}, x^\dag)
\le \frac{1}{1-\beta} \left(\frac{\d^r}{\d_*^r} \Theta(\a_*, \tilde{y}) +  \varphi\left(\d +\d_*\right)\right).
\end{equation}
In order to complete the proof, we need to estimate $\Theta(\a_*, \tilde{y})$. We will achieve this by choosing a suitable
$\hat \a\in \Delta_q$ and estimating $\Theta(\hat \a, \tilde{y})$. We fix a number $\tau\ge 6$.
If $\|F(\tilde{x}_\a)-\tilde{y}\|\le \tau \d$
for all $\a\in \Delta_q$, then we take $\hat \a = \a_0$ and obtain
\begin{equation}\label{4.14}
\Theta(\a_*, \tilde{y}) \le \Theta(\hat \a, \tilde{y}) = \frac{\|F(\tilde{x}_{\hat \a})-\tilde{y}\|^r}{\hat \a}
\le \frac{(\tau \d)^r}{\a_0}. 
\end{equation}
If there is an $\a\in \Delta_q$ such that $\|F(\tilde{x}_\a)-\tilde{y}\|>\tau \d$, we define
$\hat \a$ to be the largest number in $\Delta_q$ such that
\begin{equation*}
\|F(\tilde{x}_{q\hat \a})-\tilde{y}\| \le \tau \d <\|F(\tilde{x}_{\hat \a}) - \tilde{y}\|.
\end{equation*}
Note that the minimizing property of $\tilde{x}_\a$ implies that
$\|F(\tilde{x}_\a) -\tilde{y}\|^r \le \d^r + \a \R(x^\dag)\rightarrow \d^r$ as $\a\rightarrow 0$,
this $\hat \a$ is well-defined. Moreover, using $\tau \d <\|F(\tilde{x}_{\hat \a})-\tilde{y}\|$ and the minimizing
property of $\tilde{x}_{\hat \a}$ we can derive that
\begin{eqnarray*}
(\tau \d)^r + \hat \a \R(\tilde{x}_{\hat \a})
& \le \|F(\tilde{x}_{\hat \a}) - \tilde{y}\|^r + \hat \a \R(\tilde{x}_{\hat \a})\\
& \le \d^r + \hat \a \R(x^\dag).
\end{eqnarray*}
Since $\tau \ge 1$, we have $\R(\tilde{x}_{\hat \a}) \le \R(x^\dag)$ which implies
that $\tilde{x}_{\hat \a} \in {\mathcal M}_\rho$. Thus we may use (\ref{3.2}) with $\a =\hat \a$ to obtain
$
\tau \d \le 5 \d + \Phi^{-1} (2^r \hat \a)
$
which then implies that
\begin{equation*}
\hat \a \ge 2^{-r}\Phi((\tau-5)\d) \ge 2^{-r} \Phi(\d).
\end{equation*}
Consequently, since $\a_*$ is a global  minimizer of $\Theta$ over $\Delta_q$ and $q \hat \a\in \Delta_q$,
we can obtain
\begin{eqnarray}\label{4.15}
\Theta(\a_*,\tilde{y}) &\le \Theta(q \hat \a, \tilde{y}) =\frac{\|F(\tilde{x}_{q\hat \a})-\tilde{y}\|^r}{q \hat \a}
\le \frac{(\tau\d)^r}{q 2^{-r} \Phi(\d)}= \frac{2^{r}\tau^r }{q } \varphi(\d).
\end{eqnarray}
Combining (\ref{4.14}), (\ref{4.15}) with (\ref{3.3}) we obtain the desired estimate.
}

The a posteriori estimate in Theorem \ref{thm1} involves the quantity $\d_*$.
If $\d_*$ is about the order of $\d$, it gives convergence rates comparable
to the ones obtained in \cite{HM2012} under the Morozov's discrepancy principle. If $\d_*$ is much
larger than $\d$, only weaker convergence rates are available. If $\d_*$ is significantly smaller than $\d$,
the factor $\d/\d_*$ blows up and the approximation may diverge. Therefore, the quantity $\d_*$
provides an a posteriori check of Rule \ref{HRR}, its value should always
be monitored and the computed approximation should be discarded if $\d_*$ is presumably too small.

We can get rid of the factor $\d^r/\d_*^r$ appearing in the estimate in Theorem \ref{thm1} if the
following additional condition is stipulated on the random noise $\tilde{y}-y$. This condition also allows
to show the existence of $\a_*$ satisfying Rule \ref{HRR} and $\tilde{x}_{\a_*} \in {\mathcal M}_\rho$
which are required in Theorem \ref{thm1}.

\begin{assumption}\label{A2}
{\it There is a constant $\kappa>0$ such that
\begin{equation}\label{3.555}
\|\tilde{y}-y - v \| \ge \kappa \|\tilde{y}-y\|
\end{equation}
for any $v\in \{F(x)-y: x\in \D(F)\cap \D(\p \R)\}$.
}
\end{assumption}

Assumption \ref{A2} can be interpreted as follows. For inverse problems the forward operator $F$ usually has smoothing effect
so that $F(x)$ admits certain regularity, while the noise $\tilde{y}-y$ in general comes from randomness and hence
contains many high frequency components so that it may exhibit salient irregularity. The condition (\ref{3.555})
roughly means that subtracting any regular function of the form $F(x)-y$ from the noise can not significantly
remove the randomness.

When $\Y$ is a Hilbert space and $F$ is a bounded linear operator, it was proposed in \cite{HR1996} to use
the condition
\begin{equation*}
\|Q(\tilde{y}-y)\|\ge \sigma \|\tilde{y}-y\|
\end{equation*}
with $\sigma>0$ to prescribe the randomness of noise, where $Q$ denotes the orthogonal projection onto the orthogonal complement
of the range of $F$. This condition was weakened in \cite{JL2010} to the form: there exists $0<\sigma<1$ such that
\begin{equation}\label{4.11}
\l \tilde{y}-y, v\r \le (1-\sigma) \|\tilde{y}-y\| \|v\|
\end{equation}
for all $v\in \{ F(x)-y: x \in \D(F)\cap \D(\p \R)\}$. It is worth pointing out that (\ref{4.11}) implies Assumption \ref{A2}.
In fact, by the Cauchy-Schwarz inequality we have
\begin{eqnarray*}
\|\tilde{y}-y - v \|^2 &= \|v\|^2 + \|\tilde{y}-y\|^2 - 2 \l \tilde{y}-y, v\r \\
& \ge \|v\|^2 + \|\tilde{y}-y\|^2 - 2 (1-\sigma) \|\tilde{y}-y\|\|v\| \\
& \ge \|v\|^2 + \|\tilde{y}-y\|^2 - (1-\sigma) \left(\|\tilde{y}-y\|^2 +\|v\|^2\right)\\
& \ge \sigma \|\tilde{y}-y\|^2
\end{eqnarray*}
which shows (\ref{3.5}) with $\kappa = \sigma^{1/2}$.

\begin{corollary}\label{cor1}
Assume that $\|\tilde{y}-y\|^r \le \a_0 \R(x^\dag)$ and that $\tilde{y}-y$ satisfies Assumption \ref{A2}.
Then Rule \ref{HRR} determines a parameter $\a_*\in \Delta_q$ with the properties
\begin{equation*}
\d_*:= \|F(\tilde{x}_{\a_*})-\tilde{y}\| \ge  \kappa \d \quad \mbox{ and }
\quad \a_*\ge \frac{q\kappa^r \d^r}{(q+1) \R(x^\dag)},
\end{equation*}
where $\d = \|\tilde{y}-y\|$. If in addition $\frac{q \kappa^r+q+1}{q \kappa^r} \R(x^\dag)<\rho$, then
$\tilde{x}_{\a_*}^\d \in {\mathcal M}_\rho$  and therefore, if $x^\dag$ satisfies Assumption \ref{A1}, then
$$
D_{\xi^\dag}\R(\tilde{x}_{\a_*}, x^\dag) \le C \kappa^{-1}\left(\d^r+\varphi(\d+\d_*)\right),
$$
where $C$ is a constant depending only on $\a_0$, $q$, $r$ and $\beta$.

\end{corollary}

\proof{
From Assumption \ref{A2} it follows that
\begin{equation*}
\|F(\tilde{x}_\a) - \tilde{y} \| \ge \kappa \|y-\tilde{y}\| =\kappa \d
\end{equation*}
for all $\a>0$. This in particular shows that $\d_*\ge \kappa \d$. Furthermore, $\Theta(\a, \tilde{y}) \ge
(\kappa \d)^r/\a \rightarrow \infty$ as $\a \rightarrow 0$. This shows the existence of $\a_*$ determined
by Rule \ref{HRR}.

To derive the lower bound for $\a_*$, we first use the minimizing property of $\tilde{x}_\a$ to derive that
\begin{equation*}
\|F(\tilde{x}_\a) -\tilde{y}\|^r \le \d^r + \a \R(x^\dag), \quad \forall \a>0.
\end{equation*}
Therefore, by the definition of $\a_*$ and the lower bound on $\d_*$, we have
\begin{equation*}
\frac{(\kappa \d)^r}{\a_*} \le \Theta(\a_*, \tilde{y}) \le \Theta(\a, \tilde{y}) \le \frac{\d^r}{\a} + \R(x^\dag), \quad
\forall \a \in \Delta_q.
\end{equation*}
Now we choose $\a\in \Delta_q$ such that
\begin{equation*}
\frac{q \d^r}{\R(x^\dag)} < \a \le \frac{\d^r}{\R(x^\dag)}.
\end{equation*}
Since $\d^r \le \a_0 \R(x^\dag)$, this $\a\in \Delta_q$ is well-defined. Consequently
\begin{equation*}
\frac{(\kappa \d)^r}{\a_*} \le \left(1 + \frac{1}{q}\right) \R(x^\dag)
\end{equation*}
which implies the desired lower bound on $\a_*$. By using the minimizing property of $\tilde{x}_{\a_*}$ we then obtain
\begin{equation*}
\R(\tilde{x}_{\a_*}) \le \frac{\d^r}{\a_*} + \R(x^\dag)
\le \left(\frac{q+1}{q\kappa^r} + 1\right) \R(x^\dag)<\rho.
\end{equation*}
Thus $\tilde{x}_{\a_*} \in {\mathcal M}_\rho$. The remaining part now follows from Theorem \ref{thm1}.
}

\section{\bf Convergence}
\setcounter{equation}{0}

In Theorem \ref{thm1} and Corollary \ref{cor1} we have derived a posteriori error estimates in
terms of the Bregman distance for individually given noisy data. It is natural to ask, for a
sequence of noisy data $\{y^\d\}$ satisfying $y^\d\rightarrow y$ as $\d\rightarrow 0$, if we define
$x_\a^\d$ by
\begin{equation}\label{Tik3.1}
x_\a^\d \in \arg\min_{x\in \D(F)} \left\{ \|F(x)-y^\d\|^r + \a \R(x)\right\}
\end{equation}
and choose $\a_*:=\a_*(y^\d)$ by Rule \ref{HRR} with $\Theta(\a, \tilde{y})$ replaced by $\Theta(\a, y^\d)
:= \|F(x_\a^\d)-y^\d\|^r/\a$, i.e.
\begin{equation}\label{HR3}
\a_* \in \arg\min_{\a\in \Delta_q} \left\{\Theta(\a, y^\d):= \frac{\|F(x_\a^\d)-y^\d\|^r}{\a}\right\},
\end{equation}
is it possible to guarantee a convergence of $x_{\a_*}^\d$
to $x^\dag$ as $\d\rightarrow 0$? Bakushinskii showed in \cite{B1984} that any parameter choice rule
without using information on noise level can not guarantee a convergent regularization method for
ill-posed problems in the worst case scenario. Therefore, in order to establish a convergence result on
heuristic rules, additional conditions should be imposed on $\{y^\d\}$.
In this section we will assume
that $\{y^\d\}$ satisfies Assumption \ref{A2} uniformly in the following sense.

\begin{assumption}\label{A5}
{\it $\{y^\d\}$ is a sequence of noisy data satisfying $y^\d \rightarrow y$ as $\d\rightarrow 0$ and
there is a constant $\kappa>0$ such that
\begin{equation}\label{3.5}
\| y^\d-y - v\| \ge \kappa \|y^\d-y\|
\end{equation}
for every $y^\d$ and every $v\in \{F(x)-y: x\in \D(F)\cap \D(\p \R)\}$.
}
\end{assumption}

Under Assumption \ref{A5} we will provide four convergence results: The first one is based on the source
conditions stipulated in Assumption \ref{A1} while the other three do not depend on any source conditions.

To derive the convergence under the source conditions given in Assumption \ref{A1}, we need the
following simple fact.

\begin{lemma}\label{lem0}
Let $\Phi$ be defined by (\ref{Phi}) with $r>1$. There holds
\begin{equation*}
\lim_{t \rightarrow 0+} \frac{[\Phi^{-1}(t)]^r}{t} =0.
\end{equation*}
\end{lemma}

\proof{
Let $\gamma = [\Phi^{-1}(t)]^r/t$. Then $t= \Phi((\gamma t)^{1/r})$ which together with the definition
of $\Phi$ gives $\gamma = \varphi((\gamma t)^{1/r})$. The concavity of $\varphi$ implies that
\begin{equation*}
\varphi(t) \le C_0 t + C_1, \quad \forall t> 0
\end{equation*}
for some positive constants $C_0$ and $C_1$. Thus $\gamma \le C_0 (\gamma t)^{1/r} + C_1$ for all $t>0$.
Since $1<r<\infty$, this implies that $\gamma$ is bounded as $t \rightarrow 0$. Consequently $\gamma=
\varphi((\gamma t)^{1/r}) \rightarrow 0$ as $t\rightarrow 0$.
}

Now we are ready to give the convergence result under Assumption \ref{A1} and Assumption \ref{A5}.

\begin{theorem}\label{thm2}
Let $\{y^\d\}$ be a sequence of noisy data satisfying Assumption \ref{A5}. Let $\a_*\in \Delta_q$ be
determined by (\ref{HR3}). If $x^\dag$ satisfies Assumption \ref{A1}, then
\begin{equation*}
D_{\xi^\dag} \R(x_{\a_*}^\d, x^\dag) \rightarrow 0 \quad \mbox{ as } \d\rightarrow 0.
\end{equation*}
\end{theorem}

\proof{
We first show that $\Theta(\a_*, y^\d) \rightarrow 0$ as $\d\rightarrow 0$. By using the estimate (\ref{3.2})
and the fact that $\a_*$ is a global minimizer of $\Theta$ over $\Delta_q$, we have for all $\a\in \Delta_q$ that
\begin{eqnarray*}
\fl \Theta(\a_*,y^\d) \le \Theta(\a,y^\d) = \frac{\|F(x_\a^\d)-y^\d\|^r}{\a} \le C \left( \frac{\|y^\d-y\|^r}{\a} + \frac{[\Phi^{-1}(2^r\a)]^r}{\a}\right).
\end{eqnarray*}
Since $y^\d \rightarrow y$, we may choose $\a:=\a(\d) \in \Delta_q$ such that $\a\rightarrow 0$ and
$\|y^\d-y\|^r/\a\rightarrow 0$ as $\d\rightarrow 0$. With the help of Lemma \ref{lem0} we obtain
$\Theta(\a_*, y^\d)\rightarrow 0$ as $\d\rightarrow 0$.

In view of the facts that $\a_* \le \a_0$ and $\|F(x_{\a_*}^\d)-y^\d\|\ge \kappa \|y^\d-y\|$ we then obtain
\begin{equation*}
\fl \|F(x_{\a_*}^\d)-y^\d\|^r \le \a_0 \Theta(\a_*, y^\d) \rightarrow 0 \quad \mbox{and} \quad
\frac{\kappa^r \|y^\d-y\|^r}{\a_*} \le \Theta(\a_*, y^\d) \rightarrow 0
\end{equation*}
as $\d\rightarrow 0$. It then follows from (\ref{3.1}) that
\begin{equation*}
\fl D_{\xi^\dag}\R(x_{\a_*}^\d, x^\dag) \le C\left( \frac{\|y^\d-y\|^r}{\a_*}
+\varphi\left(\|y^\d-y\|+ \|F(x_{\a_*}^\d)-y^\d\|\right)\right) \rightarrow 0
\end{equation*}
as $\d\rightarrow 0$. This completes the proof.
}

\begin{remark}
{\rm For bounded linear operator $F$ with $\Y$ being a Hilbert space, a convergence result was proved in \cite{JL2010}
under the source condition $\xi^\dag:=F^* w \in \p \R(x^\dag)$ for some $w \in \Y$. Theorem \ref{thm2}
improves this result by showing that the convergence in fact holds under more general source conditions.
Furthermore, our proof is much simpler. This simple argument is achieved via the use of the
estimate (\ref{3.1}) which enables us to avoid the discussion on the behavior of $\a_*$ as $\d\rightarrow 0$.
}
\end{remark}

Next we will provide three convergence results without assuming any source conditions. For the first one
we need the following nonlinearity condition.

\begin{assumption}\label{A3}
{\it There is a bounded linear operator $A: \X\to \Y$ and $0\le \eta<1$ such that
\begin{equation*}
\|F(x) - F(x^\dag) - A (x-x^\dag)\| \le \eta \|F(x)-F(x^\dag)\|
\end{equation*}
for all $x\in {\mathcal M}_\rho$ with $\rho>\R(x^\dag)$.
}
\end{assumption}

Assumption \ref{A3} does not require $F$ to be Fr\'{e}chet differentiable; in case
$F$ is Fr\'{e}chet differentiable, we may take $A = F'(x^\dag)$, where $F'(x^\dag)$
denotes the Fr\'{e}chet derivative of $F$ at $x^\dag$. The condition given in Assumption \ref{A3}
is the so-called tangential cone condition which has been widely used
in the analysis of regularization methods for nonlinear inverse problems; see
\cite{HNS1995,Jin2015,JZ2014,R1999,SKHK2012} for instance.

\begin{theorem}\label{thm4}
Let $F$ satisfy Assumption \ref{A3}, let $x^\dag$ be an interior point of $\D(F)$ and let $\R$
be continuous at $x^\dag$. Let $\{y^\d\}$ be a sequence of noisy data satisfying Assumption \ref{A5} and
let $\a_*\in \Delta_q$ be determined by (\ref{HR3}). Then there exists $\xi^\dag\in \p \R(x^\dag)$
such that $D_{\xi^\dag}\R(x_{\a_*}^\d, x^\dag) \rightarrow 0$ as $\d\rightarrow 0$.
\end{theorem}

\proof{
Let $S := \{x\in {\mathcal M}_\rho: F(x) = y\}$. Clearly $x^\dag \in S$. By Assumption \ref{A3} it is
straightforward to show that
\begin{equation}\label{4.51}
S= \{ x\in {\mathcal M}_\rho: A(x-x^\dag)=0\}.
\end{equation}
According to the given conditions on $x^\dag$, we can show that the normal cone of $S$ at $x^\dag$ is
\begin{equation}\label{4.123}
N_S(x^\dag) := \{\xi \in \X^*: \l \xi, x-x^\dag\r \le 0 \ \forall x\in S\} = {\mathscr N} (A)^\perp.
\end{equation}
Indeed, since $x^\dag$ is an interior point of $\D(F)$ and $\R$ is continuous at $x^\dag$, we can find
a ball $B_\gamma(x^\dag):=\{x\in X: \|x-x^\dag\|<\gamma\}$ of radius $\gamma>0$ such that $B_\gamma(x^\dag) \subset {\mathcal M}_\rho$.
Thus, for any $x\in {\mathscr N}(A)$ we may use (\ref{4.51}) to conclude $\pm \frac{\gamma x}{\|x\|+1} +x^\dag \in S$.
Consequently $\xi\in N_S(x^\dag)$ implies that
\begin{equation*}
\pm \frac{\gamma}{\|x\|+1} \l \xi, x\r \le 0
\end{equation*}
and hence $\l \xi, x\r =0$ for all $x\in {\mathscr N}(A)$. This shows that $N_S(x^\dag) \subset {\mathscr N}(A)^\perp$
and therefore $N_S(x^\dag) = {\mathscr N}(A)^\perp$ since the opposite inclusion is obvious.

Since $X$ is reflexive, we have from (\ref{4.123}) that $N_S(x^\dag)= \overline{{\mathscr R}(A^*)}$.
Notice that
\begin{equation*}
x^\dag \in \arg\min_{x\in \X} \{ \R(x) + \iota_S(x)\},
\end{equation*}
where $\iota_S$ denotes the indicator function of $S$. Since $\R$ is continuous at $x^\dag\in S$, from
Moreau-Rockafellar theorem (\cite{Z2002}) on the sum rule of subdifferentials we have
\begin{eqnarray*}
0& \in \p(\R + \iota_S)(x^\dag) = \p \R(x^\dag) + \p \iota_S(x^\dag) = \p \R(x^\dag) + N_S(x^\dag)\\
& = \p \R(x^\dag) +  \overline{{\mathscr R}(A^*)}.
\end{eqnarray*}
Therefore there exists $\xi^\dag\in \p \R(x^\dag)$ such that $\xi^\dag \in \overline{{\mathscr R}(A^*)}$.
Thus, for any $\sigma>0$ we can find $w_\sigma\in \Y^*$ such that
\begin{equation}\label{4.31}
\|\xi^\dag - A^* w_\sigma\| \le \sigma.
\end{equation}

Now we show that $\Theta(\a_*, y^\d)\rightarrow 0$ as $\d\rightarrow 0$. To this end, we choose
$\a:=\a(\d) \in \Delta_q$ such that $\a\rightarrow 0$ and $\|y^\d-y\|^r/\a\rightarrow 0$ as
$\d\rightarrow 0$. By using the minimizing property of $x_\a^\d$ we obtain
\begin{equation}\label{4.1}
\|F(x_\a^\d) - y^\d\|^r + \a \R(x_\a^\d) \le \|y^\d-y\|^r + \a \R(x^\dag).
\end{equation}
This implies that
\begin{eqnarray*}
\|F(x_\a^\d) -y^\d \|^r + \a D_{\xi^\dag} \R(x_\a^\d, x^\dag)
& \le \|y^\d-y\|^r - \a \l \xi^\dag, x_\a^\d-x^\dag\r \\
& = \|y^\d-y\|^r - \a \l\xi^\dag -A^*w_\sigma, x_\a^\d-x^\dag\r\\
& \quad \, - \a \l w_\sigma, A(x_\a^\d -x^\dag)\r.
\end{eqnarray*}
Consequently, by virtue of (\ref{4.31}) and Assumption \ref{A3} we have
\begin{eqnarray*}
\|F(x_\a^\d) -y^\d \|^r + \a D_{\xi^\dag} \R(x_\a^\d, x^\dag)
& \le  (1+\eta) \a \|w_\sigma\| \|F(x_\a^\d) -y\| \\
& \quad \, + \|y^\d-y\|^r + \a \sigma \|x_\a^\d-x^\dag\|.
\end{eqnarray*}
According to (\ref{4.1}) we have $\|F(x_\a^\d)-y\|\rightarrow 0$ and $\R(x_{\a}^\d) \le \|y^\d-y\|^r/\a+\R(x^\dag)
\rightarrow \R(x^\dag)$ as $\d\rightarrow 0$. Thus, by the coercivity of the function $x\to \|F(x)-y\|^r + \R(x)$
we can conclude that $\|x_{\a}^\d\|$ is bounded and hence $\|x_{\a}^\d-x^\dag\| \le C_0$ for some
constant $C_0$ independent of $\d$. Consequently
\begin{equation*}
\|F(x_{\a}^\d) -y^\d \|^r \le \|y^\d-y\|^r + C_0 \a \sigma + (1+\eta) \a \|w_\sigma\| \|F(x_\a^\d)-y\|.
\end{equation*}
This implies that
\begin{eqnarray*}
\limsup_{\d\rightarrow 0} \Theta(\a, y^\d)
& \le \lim_{\d\rightarrow 0} \left(\frac{\|y^\d-y\|^r}{\a} + C_0 \sigma + (1+\eta) \|w_\sigma\| \|F(x_\a^\d)-y\|\right)\\
& = C_0 \sigma
\end{eqnarray*}
Because $\a_*$ is a global minimizer of $\Theta$ over $\Delta_q$, we can obtain
\begin{equation*}
\limsup_{\d\rightarrow 0} \Theta(\a_*, y^\d) \le \limsup_{\d\rightarrow 0} \Theta(\a, y^\d)
\le C_0 \sigma.
\end{equation*}
Since $\sigma>0$ can be arbitrarily small and $\Theta(\a_*, y^\d)$ is nonnegative, we must have
$\lim_{\d\rightarrow 0} \Theta(\a_*, y^\d)=0$. This together with the facts that $\a_*\le \a_0$ and
$\|F(x_{\a_*}^\d) -y^\d\|\ge \kappa \|y^\d-y\|$ from Corollary \ref{cor1} shows that
\begin{equation}\label{4.2}
\|F(x_{\a_*}^\d) - y^\d \| \rightarrow 0 \quad \mbox{and} \quad
\frac{\|y^\d-y\|^r}{\a_*} \rightarrow 0 \quad \mbox{ as } \d\rightarrow 0.
\end{equation}

Finally we prove $D_{\xi^\dag} \R(x_{\a_*}^\d, x^\dag) \rightarrow 0$ as $\d\rightarrow 0$.
Because of  (\ref{4.2}), we may use (\ref{4.1}) to show the boundedness of $\{\R(x_{\a_*}^\d)\}$
and $\{F(x_{\a_*}^\d)\}$ which together with the coercivity of $x\to \|F(x)-y^\d\|^r + \R(x)$ shows
the boundedness of $\{x_{\a_*}^\d \}$.
By taking a subsequence if necessary, we can conclude that $x_{\a_*}^\d \rightharpoonup \hat x$ for
some $\hat x \in \X$ as $\d\rightarrow 0$. In view of (\ref{4.1}), (\ref{4.2}) and the lower semi-continuity
of norms and $\R$ we can derive that
\begin{equation*}
0\le \|F(\hat x) -y\| \le \lim_{\d\rightarrow 0} \|F(x_{\a_*}^\d) -y^\d\| =0
\end{equation*}
and
\begin{eqnarray*}
\R(\hat x) & \le \liminf_{\d\rightarrow 0} \R(x_{\a_*}^\d)
\le \limsup_{\d\rightarrow 0} \R(x_{\a_*}^\d) \\
& \le \lim_{\d\rightarrow 0} \left(\frac{\|y^\d-y\|^r}{\a_*} + \R(x^\dag)\right)
= \R(x^\dag).
\end{eqnarray*}
Thus $F(\hat x) = y$. Since $x^\dag$ is an $\R$-minimizing solution of $F(x) = y$ in ${\mathcal M}_\rho$, we have
$\R(\hat x) =  \R(x^\dag)$ and hence
\begin{equation*}
\lim_{\d\rightarrow 0} \R(x_{\a_*}^\d) = \R(x^\dag).
\end{equation*}
This together with the fact $x_{\a_*}^\d \rightharpoonup \hat x$ shows that
\begin{equation*}
\lim_{\d\rightarrow 0} D_{\xi^\dag}\R(x_{\a_*}^\d, x^\dag)
= \lim_{\d\rightarrow 0}  \l \xi^\dag, x^\dag-x_{\a_*}^\d\r =\l \xi^\dag, x^\dag-\hat x\r.
\end{equation*}
Since $\xi^\dag \in {\mathscr N}(A)^\perp$ and $\hat x-x^\dag \in {\mathscr N}(A)$ we must have
$\lim_{\d\rightarrow 0} D_{\xi^\dag}\R(x_{\a_*}^\d, x^\dag) =0$. The proof is thus complete.
}

\begin{remark}\label{rk3.2}
{\rm In Theorem \ref{thm4} we obtain the convergence of $x_{\a_*}^\d$ to $x^\dag$ in the Bregman distance. This
does not imply the convergence in norm in general. However, from the proof of Theorem \ref{thm4}
it is easily seen that, if $x^\dag$ is the unique $\R$-minimizing solution of $F(x)=y$ in ${\mathcal M}_\rho$, we have
actually shown that $x_{\a_*}^\d \rightharpoonup x^\dag$ and $\R(x_{\a_*}^\d) \rightarrow \R(x^\dag)$
as $\d\rightarrow 0$. Consequently $\|x_{\a_*}^\d-x^\dag\| \rightarrow 0$ as $\d\rightarrow 0$ as long as
$\R$ admits the {\it Kadec property} in the sense that any sequence $\{x_n\}$ satisfying $x_n\rightharpoonup \hat x$
and $\R(x_n) \rightarrow \R(\hat x)<\infty$ must have $\|x_n-\hat x\|\rightarrow 0$ as $n\rightarrow \infty$.
}
\end{remark}

\begin{remark}
{\rm Due to the relation (\ref{4.51}), one can show under Assumption \ref{A3} that $x^\dag$ is the unique
$\R$-minimizing solution of $F(x) =y$ in ${\mathcal M}_\rho$ if $\R$ is strictly convex on ${\mathscr N}(A)$.
}
\end{remark}

The convergence result given in Theorem \ref{thm4} requires $\R$ to be continuous at least at one point in $S$.
This condition is already very weak. However, there are important situations for which $\R$ is nowhere continuous.
The typical examples are the $\ell^1$-norm $\R(x) = \|x\|_{\ell^1}$ in the sequence space $\ell^2$ and the
total variation
\begin{equation*}
\R(x) = \int_\Omega |\nabla x| := \sup\left\{\int_\Omega x \mbox{div} f d\mu: f\in C_0^1(\Omega, {\mathbb R}^d)
\mbox{ and } \|f\|_{L^\infty} \le 1\right\}
\end{equation*}
in the function space $L^2(\Omega)$, where $\Omega\subset {\mathbb R}^d$ is a bounded Lipschitz domain.

Our next two results provide convergence criteria without assuming continuity on $\R$.
The first one requires the linear operator $A$ in Assumption \ref{A3} to be injective. We will make use of
the $\ep$-subdifferential calculus. For any $\ep>0$ the set
\begin{equation*}
\p_\ep\R(x) := \{\xi\in X^*: \R(\bar x) \ge \R(x) +\l \xi, \bar x-x\r - \ep \ \mbox{ for all } \bar x \in X\}
\end{equation*}
is called the $\ep$-subdifferential of $\R$ at $x$. We have (see \cite[Theorem 2.4.4]{Z2002}).

\begin{lemma}\label{subdiff}
If $\R: X \to (-\infty, \infty]$ is a proper, lower semi-continuous, convex function, then $\p_\ep \R(x) \ne \emptyset$
for any $x\in \D(\R)$ and $\ep>0$.
\end{lemma}

\begin{theorem}\label{thm8}
Let $F$ satisfy Assumption \ref{A3} with $A$ injective and let $x^\dag$ be the unique $\R$-minimizing
solution of (\ref{1}) in ${\mathcal M}_\rho$. Let $\{y^\d\}$ be a sequence of noisy data
satisfying Assumption \ref{A5} and let $\a_*\in \Delta_q$ be determined by (\ref{HR3}). Then
\begin{equation*}
x_{\a_*}^\d \rightharpoonup x^\dag, \quad \R(x_{\a_*}^\d)\rightarrow \R(x^\dag) \quad \mbox{and} \quad
F(x_{\a_*}^\d) \rightarrow y
\end{equation*}
as $\d\rightarrow 0$. If, in addition,  $\R$ admits the Kadec property, then $x_{\a_*}^\d \rightarrow x^\dag$
as $\d\rightarrow 0$.
\end{theorem}

\proof{
According to the proof of Theorem \ref{thm4} and Remark \ref{rk3.2}, it suffices to show that $\Theta(\a_*, y^\d)
\rightarrow 0$ as $\d\rightarrow 0$.  Since $y^\d\rightarrow y$, we may choose $\a:=\a(\d)\in \Delta_q$
such that $\a \rightarrow 0$ and $\|y^\d-y\|^r/\a\rightarrow 0$ as $\d \rightarrow 0$. For any $\ep>0$, we may use
Lemma \ref{subdiff} to find an element $\xi_\ep\in \p_\ep \R(x^\dag)$. By making use of (\ref{4.1}) we have
\begin{equation*}
\|F(x_\a^\d)-y^\d\|^r + \a D_{\xi_\ep}^\ep \R(x_\a^\d, x^\dag)
\le \|y^\d-y\|^r + \a \ep - \a \l \xi_\ep, x_\a^\d-x^\dag\r,
\end{equation*}
where
\begin{equation*}
D_{\xi_\ep}^\ep \R(x_\a^\d, x^\dag) := \R(x_\a^\d) -\R(x^\dag) -\l \xi_\ep, x_\a^\d-x^\dag\r +\ep
\end{equation*}
which is nonnegative. Since $A$ is injective, we have $X^* = {\mathscr N}(A)^\perp = \overline{{\mathscr R}(A^*)}$.
Thus, for any $\sigma>0$ we can find $w_\sigma\in Y^*$ such that $\|\xi_\ep - A^* w_\sigma\| \le \sigma$. Therefore
\begin{eqnarray*}
\fl \quad \|F(x_\a^\d)-y^\d\|^r
&\le \|y^\d-y\|^r + \a \ep - \a \l \xi_\ep-A^* w_\sigma, x_\a^\d-x^\dag\r -\a \l w_\sigma, A(x_\a^\d-x^\dag)\r\\
&\le \|y^\d-y\|^r + \a \ep + \a\sigma \|x_\a^\d-x^\dag\| +\a \| w_\sigma\| \| A(x_\a^\d-x^\dag)\|.
\end{eqnarray*}
As in the proof of Theorem \ref{thm4} we can find a universal constant $C$ such that
\begin{eqnarray*}
\|F(x_\a^\d)-y^\d\|^r
&\le \|y^\d-y\|^r + \a \ep + C \a\sigma + C \a \| w_\sigma\| \| F(x_\a^\d)-y\|.
\end{eqnarray*}
Using the choice of $\a$ we have $\|F(x_\a^\d)-y\| \rightarrow 0$ and thus
\begin{equation*}
\limsup_{\d\rightarrow 0} \Theta(\a_*, y^\d) \le \limsup_{\d\rightarrow 0} \Theta(\a, y^\d) \le \ep + C \sigma.
\end{equation*}
Since $\sigma>0$ and $\ep>0$ can be arbitrarily small, we therefore obtain $\Theta(\a_*, y^\d) \rightarrow 0$ as
$\d\rightarrow 0$.
}

Finally we give a convergence result which use neither the continuity of $\R$ nor the injectivity of
the linearized operator of $F$ at $x^\dag$. However, we need to
restrict to the situation that $\Y$ is a Hilbert space and $r=2$ in the formulation of (\ref{Tik3.1}).
We will use $x_\a$ to denote a minimizer of (\ref{Tik3.1}) with $y^\d$ replaced by $y$, i.e.
\begin{equation*}
x_\a \in \arg\min_{x\in \D(F)} \left\{\|F(x) - y\|^2 +\a \R(x)\right\}.
\end{equation*}
We will assume that $F$ is Fr\'{e}chet differentiable and satisfies the following nonlinearity condition.

\begin{assumption}\label{A4}
{\it There exist $\rho >\R(x^\dag)$ and $\kappa \ge 0$ such that
\begin{equation*}
\|F(\bar x) -F(x) - F'(x) (\bar x-x) \| \le \kappa [D_{\xi}\R(\bar x, x)]^{1/2} \|F(\bar x) - F(x) \|
\end{equation*}
for all $\bar x, x\in {\mathcal M}_\rho$ and $\xi \in \p \R(x)$.
}
\end{assumption}

Assumption \ref{A4} has been used in the work of regularization theory for nonlinear ill-posed inverse problems
in Banach spaces, see for instance \cite{HY2010,JZ2014} and the references therein. When $\R$ is $2$-convex
in the sense that there is a constant $C_0$ such that
\begin{equation*}
\|\bar x- x\|^2 \le C_0 D_{\xi} \R(\bar x, x)
\end{equation*}
for all $\bar x, x \in \D(\R)$ and $\xi \in \p \R(x)$, Assumption \ref{A4} holds if there is a constant
$\kappa\ge 0$ such that
\begin{equation*}
\|F(\bar x)-F(x) -F'(x) (\bar x-x)\| \le \kappa \|\bar x-x\| \|F(\bar x)-F(x)\|
\end{equation*}
for all $\bar x, x\in {\mathcal M}_\rho$, which is a slightly strengthened version of Assumption \ref{A3}.

\begin{lemma}\label{lem5}
Let $F$ be Fr\'{e}chet differentiable and satisfy Assumption \ref{A4}. Assume that $x_\a^\d, x_\a \in {\mathcal M}_\rho$
and ${\mathcal M}_\rho$ is contained in the interior of $\D(F)$. Then
\begin{equation*}
\fl \|F(x_\a^\d) -y^\d + y -F(x_\a)\|^2 + 2 \a \left(1- \frac{4\kappa^2 \|F(x_\a)-y\|^2}{\a}\right) D_{\xi_\a}\R(x_\a^\d, x_\a)
\le 3 \|y^\d-y\|^2,
\end{equation*}
where $\xi_\a := \frac{2}{\a} F'(x_\a)^* (y-F(x_\a)) \in \p \R(x_\a)$.
\end{lemma}

\proof{
Since $x_\a\in {\mathcal M}_\rho$ is an interior point of $\D(F)$, the first order optimality condition shows that
$\xi_\a \in \p \R(x_\a)$. By the minimizing property of $x_\a^\d$ we have
\begin{equation*}
\|F(x_\a^\d) - y^\d\|^2 +\a \R(x_\a^\d) \le \|F(x_\a)-y^\d\|^2 +\a \R(x_\a)
\end{equation*}
which, after rearrangement, gives
\begin{eqnarray*}
\fl \|F(x_\a^\d) -y^\d + y -F(x_\a)\|^2 + \a D_{\xi_\a} \R(x_\a^\d, x_\a) \\
 \le \|y-y^\d\|^2 + 2 \l y-F(x_\a), F(x_\a^\d)- F(x_\a) \r -  \a \l \xi_\a, x_\a^\d-x_\a\r\\
 = \|y^\d-y\|^2 + 2 \l y-F(x_\a), F(x_\a^\d)-F(x_\a) - F'(x_\a) (x_\a^\d-x_\a)\r.
\end{eqnarray*}
By making use of Assumption \ref{A4} and the Cauchy-Schwarz inequality we obtain
\begin{eqnarray*}
\fl  \|F(x_\a^\d)-y^\d + y-F(x_\a)\|^2 +\a D_{\xi_\a} \R(x_\a^\d, x_\a) \\
 \le \|y^\d-y\|^2 + 2 \kappa \|y-F(x_\a)\| [D_{\xi_\a}\R(x_\a^\d, x_\a)]^{1/2} \|F(x_\a^\d) - F(x_\a)\|\\
 \le \|y^\d-y\|^2 + \frac{1}{4} \|F(x_\a^\d)-F(x_\a)\|^2 + 4 \kappa^2 \|y-F(x_\a)\|^2 D_{\xi_\a} \R(x_\a^\d, x_\a)\\
 \le \frac{3}{2} \|y^\d-y\|^2 + \frac{1}{2} \|F(x_\a^\d)-y^\d+y-F(x_\a)\|^2 \\
 \quad \, + 4 \kappa^2 \|y-F(x_\a)\|^2 D_{\xi_\a} \R(x_\a^\d, x_\a).
\end{eqnarray*}
This shows the desired inequality.
}

\begin{lemma}\label{lem6}
There holds $\|F(x_\a)-y\|^2/\a \rightarrow 0$ as $\a\rightarrow 0$.
\end{lemma}

\proof{
Let $x^\dag$ denote an $\R$-minimizing solution of (\ref{1}). By repeating the argument in the last part of the
proof of Theorem \ref{thm4}, we can obtain $\R(x_\a) \rightarrow \R(x^\dag)$ as $\a \rightarrow 0$. This fact,
together with the inequality
\begin{equation}\label{4.8}
\|F(x_\a)-y\|^2 +\a \R(x_\a) \le \a \R(x^\dag)
\end{equation}
obtained from the minimizing property of $x_\a$, shows that $\|F(x_\a)-y\|^2/\a \rightarrow 0$ as $\a\rightarrow 0$.
}

\begin{theorem}\label{thm5}
Let $\Y$ be a Hilbert space and $r=2$. Assume that $F$ is Fr\'{e}chet differentiable and satisfies Assumption
\ref{A4}. Assume also that $x^\dag$ is the unique $\R$-minimizing solution of (\ref{1}) in ${\mathcal M}_\rho$
and ${\mathcal M}_\rho$ is contained in the interior of $\D(F)$. Let $\{y^\d\}$ be a sequence of noisy data satisfying
Assumption \ref{A5}. Then for the parameter $\a_*\in \Delta_q$ determined by (\ref{HR3}) there hold
\begin{equation*}
x_{\a_*}^\d \rightharpoonup x^\dag, \quad \R(x_{\a_*}^\d)\rightarrow \R(x^\dag) \quad \mbox{and} \quad
F(x_{\a_*}^\d) \rightarrow y
\end{equation*}
as $\d\rightarrow 0$. If, in addition,  $\R$ admits the Kadec property, then $x_{\a_*}^\d \rightarrow x^\dag$
as $\d\rightarrow 0$.
\end{theorem}

\proof{
We choose $\a:=\a(\d) \in \Delta_q$ such that $\a\rightarrow 0$ and $\|y^\d-y\|^2/\a\rightarrow 0$ as $\d\rightarrow 0$.
The minimizing property of $x_\a^\d$ and $x_\a$ shows that
\begin{equation*}
\R(x_\a^\d) \le \frac{\|y^\d-y\|^2}{\a} + \R(x^\dag) \quad \mbox{and} \quad \R(x_\a) \le \R(x^\dag)
\end{equation*}
which then imply that $x_\a^\d, x_\a \in {\mathcal M}_\rho$ for small $\d>0$. Furthermore, from Lemma \ref{lem6} it follows
that $\|F(x_{\a})-y\|^2/\a\rightarrow 0$ as $\d\rightarrow 0$. Thus we may use Lemma \ref{lem5} to derive that
\begin{equation*}
\|F(x_{\a}^\d)-y^\d+y -F(x_{\a})\|^2 \le 3 \|y^\d-y\|^2
\end{equation*}
for small $\d>0$. Consequently
\begin{eqnarray*}
\Theta(\a, y^\d) &\le \frac{2\|F(x_{\a}^\d)-y^\d+y -F(x_{\a})\|^2}{\a} +\frac{2\|F(x_{\a})-y\|^2}{\a}\\
& \le \frac{6 \|y^\d-y\|^2}{\a} + \frac{2\|F(x_{\a})-y\|^2}{\a} \rightarrow 0
\end{eqnarray*}
as $\d\rightarrow 0$. Since $\a_*$ is a global minimizer of $\Theta$ over $\Delta_q$, we obtain
$\Theta(\a_*, y^\d) \le \Theta(\a, y^\d) \rightarrow 0$ as $\d\rightarrow 0$ which, together with the facts
$\a_*\le \a_0$ and $\|F(x_{\a_*}^\d)-y^\d\| \ge \kappa \|y^\d-y\|$, shows that
\begin{equation*}
\|F(x_{\a_*}^\d)-y^\d\| \rightarrow 0 \quad \mbox{ and } \quad \frac{\|y^\d-y\|^2}{\a_*} \rightarrow 0
\quad \mbox{ as } \d\rightarrow 0.
\end{equation*}
Now we can repeat the argument in the last part of the proof of Theorem \ref{thm4} to complete the proof.
}

\begin{remark}
{\rm Under the conditions in Theorem \ref{thm5} with $\X=\ell^2$ and $\R(x) = \|x\|_{\ell^1}$,
we have $x_{\a_*}^\d\rightharpoonup x^\dag$ in $\ell^2$ and $\|x_{\a_*}^\d\|_{\ell^1} \rightarrow \|x^\dag\|_{\ell^1}$
as $\d\rightarrow 0$. In view of the Kadec property of $\R$ shown in \cite[Lemma 4.3]{DDD2004} we obtain
$\|x_{\a_*}^\d - x^\dag\|_{\ell^2} \rightarrow 0$ as $\d\rightarrow 0$. By using  \cite[Lemma 2]{GHS2008} we can
even obtain the stronger result $\|x_{\a_*}^\d -x^\dag\|_{\ell^1} \rightarrow 0$ as $\d\rightarrow 0$.
}
\end{remark}

\section{\bf Numerical results}
\setcounter{equation}{0}

In this section we will provide numerical examples to test the performance of the variational regularization
(\ref{Tik}) when the regularization parameter $\a>0$ is chosen by Rule \ref{HRR}. In the following computation
all the minimization problems are solved by a gradient descent method.

\begin{figure}[ht!]
     \begin{center}
        {
           \includegraphics[width = 1\textwidth, height = 0.7\textwidth]{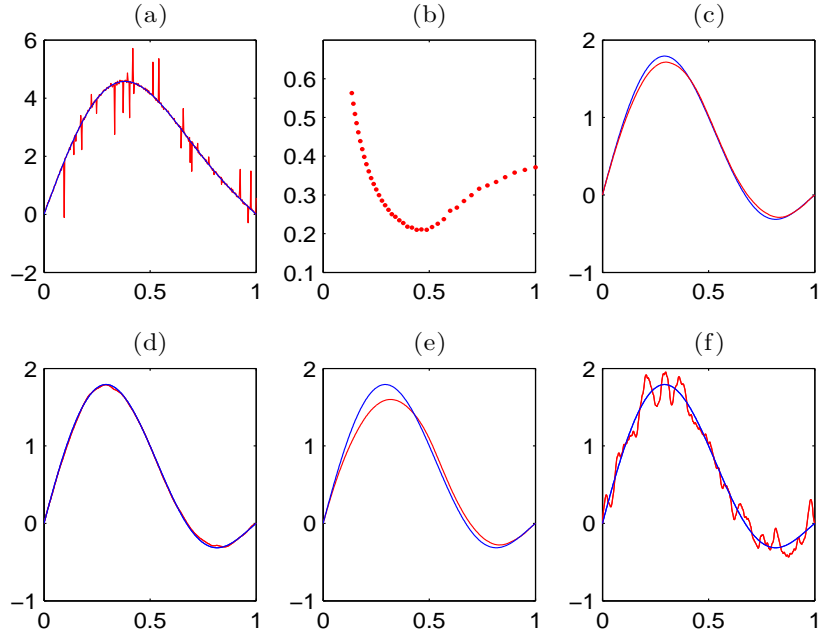}
        }
    \end{center}
    \vskip -0.3cm
    \caption{(a) noisy data with outliers; (b) $\Theta(\a, \tilde{y})$ versus $\a$; (c) reconstruction result by Rule \ref{HRR};
    (d), (e) and (f) are reconstruction results by the discrepancy principle with $\tau = 1.01, 1.615$ and $0.996$   }%
   \label{fig1}
\end{figure}

\begin{example}
{\rm We consider the linear integral equation of the form
\begin{equation*}
(F x)(s) := \int_0^1 k(s, t) x(t) dt = y(s)   \quad \mbox{ on } [0,1],
\end{equation*}
where
\begin{equation*}
k(s, t) = \left\{\begin{array}{lll}
40 s(1-t) & \mbox{ if } s\le t,\\
40 t(1-s) & \mbox{ if } t<s.
\end{array}\right.
\end{equation*}
Assume that the sought solution is
$
x^\dag (t) = 4t(1-t)+ \sin(2\pi t)
$
and the exact data $y := F x^\dag$ is corrupted by impulsive noise so that we have the noisy data
$\tilde{y}$ as shown in Figure \ref{fig1} (a). In order to use $\tilde{y}$ to reconstruct $x^\dag$ we use
the variational regularization (\ref{Tik}) with $\X= L^2[0,1]$, $\Y = L^r[0,1]$ with $r=1.01$ and
$\R(x) = \|x\|_{L^2}^2$. We choose the regularization parameter $\a$ by Rule \ref{HRR}
with $\a_0 =1$ and $ q = 0.95$. The relation between $\Theta(\a, \tilde{y})$ and $\a$ is plotted in
Figure \ref{fig1} (b) and the reconstruction result is shown in Figure \ref{fig1} (c). As comparison
we also consider the choice of the regularization parameter by the discrepancy principle which chooses
$\a$ to be the largest number in $\Delta_q$ satisfying $\|F \tilde{x}_\a -\tilde{y}\|_{L^r}\le \tau \d$, where
$\d=\|\tilde{y}-y\|_{L^r}$ is the noise level and $\tau>0$ is a given number. In Figure \ref{fig1} (d), (e) and (f)
we plot the reconstruction results for $\tau=1.01$, $\tau = 1.615$ and $\tau = 0.996$ respectively;
these three choices of $\tau$ correspond the proper estimation, overestimation and underestimation
of the noise level.  From Figure \ref{fig1} it can be seen that  Rule \ref{HRR} gives satisfactory
reconstruction result although no information on noise level is used. The discrepancy principle can give better
result if accurate information on noise level is used; however, it can give much worse result if the
noise level is overestimated or underestimated. In particular, it can produce very oscillatory result
if an underestimated noise level is used.
}
\end{example}

\begin{example}
{\rm We next consider the estimation of the coefficient $c$ in the boundary value problem
\begin{equation}\label{6.1}
\left\{\begin{array}{lll}
-u''+ c u = f \quad \mbox{ in } (0, 1)\\
u(0) = g_0, \quad u(1) = g_1
\end{array}\right.
\end{equation}
from the measurement of the state variable $u$, where $g_0$, $g_1$ and $f\in H^{-1}[0, 1]$ are given.
It is well known that (\ref{6.1}) has a unique solution $u := u(c) \in H^1[0, 1]$ for each $c$ in the domain
\begin{equation*}
\D := \{c \in L^2[0, 1]:  \|c - \hat c\|_{L^2} \le \gamma \mbox{ for some } \hat c \ge 0 \mbox{ a.e.}\}
\end{equation*}
with some $\gamma > 0$. We consider the problem of identifying $c\in L^2[0,1]$ from an $L^2[0,1]$-measurement
$\tilde{u}$ of $u$. By taking $\X = \Y = L^2[0, 1]$, this inverse problem reduces to solving
(\ref{1}) with the nonlinear operator $F : \D \subset  L^2[0, 1] \to L^2[0, 1]$ defined as $F(c) := u(c)$.
It is easy to show that $F$ is Fr\'{e}chet differentiable, and the Fr\'{e}chet derivative and its adjoint
are given by ´
\begin{equation*}
F'(c)h = −A(c)^{-1}(hu(c)), \qquad F'(c)^* w = -u(c)A(c)^{-1} w,
\end{equation*}
where $A(c) : H^2 \cap H_0^1 \to L^2$ is defined by $A(c)u = -u'' +cu$. We will reconstruct the sought
coefficient $c^\dag$ using (\ref{Tik}) with $r=2$ and various choices of $\R$ according to the available
a priori information on $c^\dag$. In all the examples we use $g_0=1$, $g_1 = 6$ and $f(t) = 100 \exp(-10(t-0.5)^2)$.
In our numerical computation, all differential equations are solved approximately by the finite
difference method by dividing the interval $[0, 1]$ into $N = 400$ subintervals of equal length.

In Figure \ref{fig2} we report the computational result of Rule \ref{HRR} with $\a_0=0.005$ and $q=0.8$ when the
sought solution is $c^\dag = \sin(\pi t) + \sin(4\pi t) + 2t^3(1-t)+t$ which is smooth. Assuming that the data is
corrupted by Gaussian noise with $\|\tilde{u}-u\|_{L^2[0,1]} = \d$ and $\d=0.0025$, we reconstruct the sought solution
by using (\ref{Tik}) with two choices of $\R$, i.e. $\R_1(c) = \|c\|_{L^2}^2$ and $\R_2(c) = \|c-c_0\|_{L^2}^2$ with $c_0(t) =t$.
According to \cite{EKN1989} we have $\p \R_2(c^\dag) \cap {\mathscr R}(F'(c^\dag)^*)\ne \emptyset$ and
this does not hold for $\R_1$. In view of the a posteriori error estimates given in Corollary \ref{cor1}
we expect that (\ref{Tik}) with $\R=\R_2$ can give better reconstruction result than $\R=\R_1$.
This is confirmed by the plots in Figure \ref{fig2} which also shows that reasonable reconstruction results
can be obtained even though the source condition (\ref{sc}) does not hold.

\begin{figure}[ht!]
     \begin{center}
        {
           \includegraphics[width = 1\textwidth, height = 0.7\textwidth]{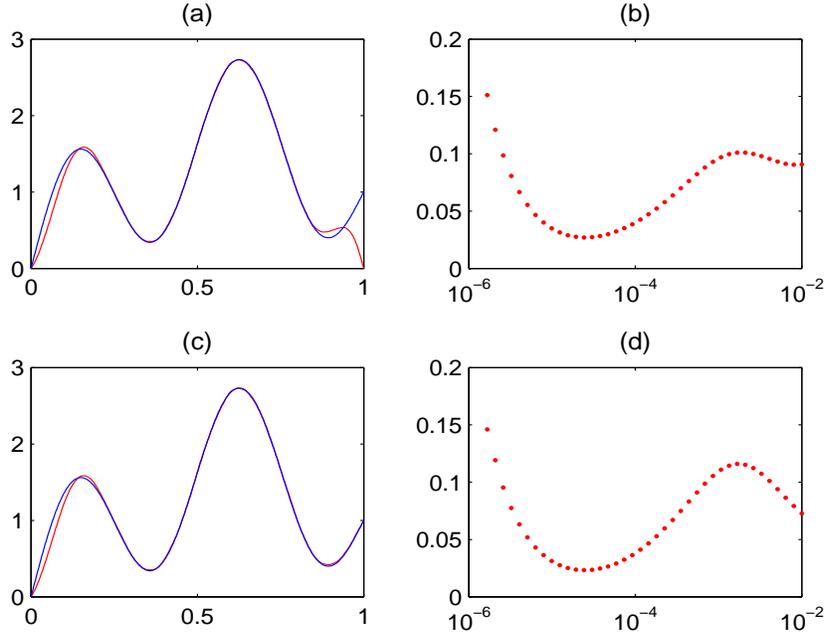}
        }
    \end{center}
    \vskip -0.3cm
    \caption{(a) reconstruction result using $\R_1$; (b) $\Theta(\a, \tilde{u})$ versus $\a$ corresponding to
    $\R =\R_1$; (c) reconstruction result using $\R_2$;
    (d) $\Theta(\a, \tilde{u})$ versus $\a$ corresponding to $\R =\R_2$.  }%
   \label{fig2}
\end{figure}

Next we consider the reconstruction performance when the sought coefficient is piecewise constant. The sought solution is
plotted in Figure \ref{fig3} (a). Assuming that the data is corrupted by Gaussian noise with $\|\tilde{u}-u\|_{L^2[0,1]} = \d$
and $\d=0.001$, we reconstruct the sought solution by using (\ref{Tik}) with
$
\R(c) = \int_{[0,1]} |D c| + \mu \|c\|_{L^2[0,1]}^2
$
and $\mu= 0.001$ in which the regularization parameter is chosen by Rule \ref{HRR} with $\a_0 = 0.001$ and $q = 0.8$.
Figure \ref{fig3} (a) and (b) plot the reconstruction result  and the curve of $\Theta(\a, \tilde{u})$ versus $\a$ respectively.
The reconstruction coincides with the sought solution very well.
}
\end{example}

\begin{figure}[ht!]
     \begin{center}
        {
        \includegraphics[width = 1\textwidth]{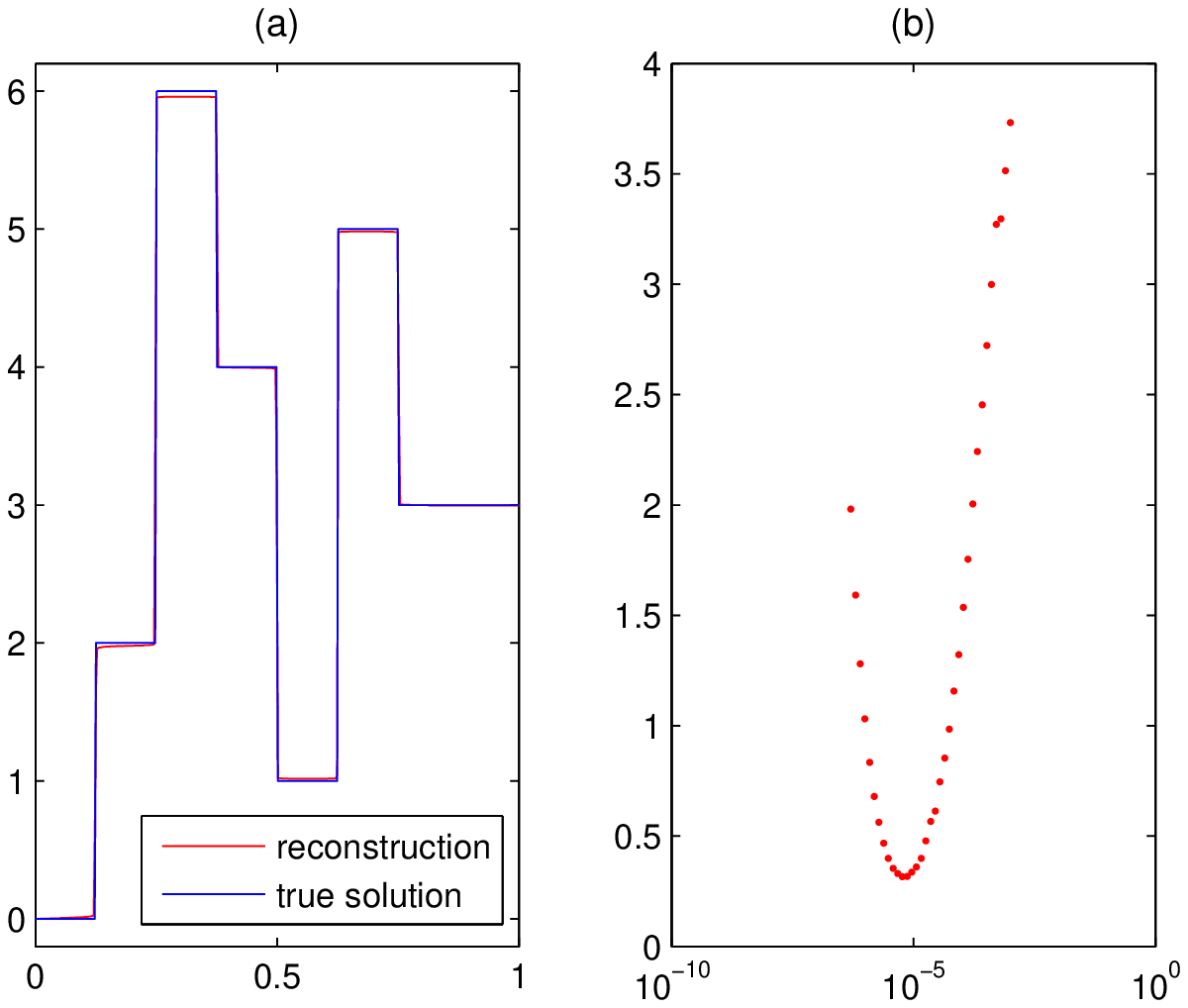}
        }
    \end{center}
    \vskip -0.3cm
    \caption{(a) reconstruction result; (b) $\Theta(\a, \tilde{u})$ versus $\a$.  }%
   \label{fig3}
\end{figure}

\section*{\bf Acknowledgement}

This work is partially supported by the Discovery Project grant DP150102345 of Australian
Research Council.

%


\section*{References}

\begin{thebibliography}{999}
%
%
\bibitem{B1984} A. B. Bakushinskii, {\it Remarks on choosing a regularization parameter using the quasi–optimality
and ratio criterion}, U.S.S.R. Comp. Maths. Math. Phys. 24 (1984), No. 4, 181--182.

\bibitem{DDD2004} I. Daubechies, M. Defrise and C. De Mol,  {\it An iterative thresholding algorithm for linear
inverse problems with a sparsity constraint},  Comm. Pure Appl. Math. 57 (2004), 1413--1457.

\bibitem{EHN1996} H. W. Engl, M. Hanke and A. Neubauer, {\it Regularization of Inverse Problems},
Vol. 375, Springer Science \& Business Media, 1996.

\bibitem{EKN1989} H. W. Engl, K. Kunisch and A. Neubauer, {\it Convergence rates for Tikhonov regularisation of
nonlinear ill-posed problems}, Inverse Problems 5 (1989), no. 4, 523--540.

\bibitem{GHS2008} M. Grasmair, M. Haltmeier and O. Scherzer, {\it Sparse regularization with $\ell^q$ penalty term},
Inverse Problems 24 (2008), no. 5, 055020, 13 pp.

\bibitem{HNS1995} M. Hanke, A. Neubauer and O. Scherzer, {\it A convergence analysis of the Landweber iteration
for nonlinear ill-posed problems}, Numer. Math. 72 (1995), 21--37.

\bibitem{HR1996} M. Hanke and T. Raus, {\it A general heuristic for choosing the regularization parameter in
ill-posed problems}, SIAM J. Sci. Comput., 17 (1996), no. 4, 956--972.

\bibitem{HL1993} P. C. Hansen and D. P. OʾLeary, {\it The use of the L-curve in the regularization of discrete
ill-posed problems}, SIAM J. Sci. Comput., 14 (1993), 1487--1503.

\bibitem{HM2012} B. Hofmann and P. Math\'{e}, {\it Parameter choice in Banach space regularization under
variational inequalities}, Inverse Problems, 28 (2012), no. 10, 104006, 17 pp.

\bibitem{HKPS2007} B. Hofmann, B. Kaltenbacher, C. Pöschl and O. Scherzer, {\it A convergence
rates result for Tikhonov regularization in Banach spaces with non-smooth operators}, Inverse Problems,
23 2007), 987--1010.

\bibitem{HY2010} B. Hofmann and M. Yamamoto, {\it On the interplay of source conditions and variational inequalities
for nonlinear ill-posed problems}, Appl. Anal., 89 (2010), no. 11, 1705--1727.

\bibitem{IJ2011} K. Ito and B. Jin, {\it A new approach to nonlinear constrained Tikhonov regularization},
Inverse Problems 27 (2011), no. 10, 105005, 23 pp.

\bibitem{JL2010} B. Jin and D. A. Lorenz, {\it Heuristic parameter-choice rules for convex variational
regularization based on error estimates}, SIAM J. Numer. Anal., 48 (2010), no. 3, 1208--1229.

\bibitem{Jin1999} Q. Jin, {\it Applications of the modified discrepancy principle to Tikhonov regularization
of nonlinear ill-posed problems}, SIAM J. Numer. Anal. 36 (1999), no. 2, 475--490.

\bibitem{Jin2015} Q. Jin, {\it Inexact Newton-Landweber iteration in Banach spaces with non-smooth convex penalty terms},
SIAM J. Numer. Anal., 53 (2015), no. 5, 2389--2413.

\bibitem{JH1999} Q. Jin and Z. Hou, {\it On an a posteriori parameter choice strategy for Tikhonov
regularization of nonlinear ill-posed problems}, Numer. Math., 83 (1999), no. 1, 139--159.

\bibitem{JZ2014} Q. Jin and M. Zhong, {\it Nonstationary iterated Tikhonov regularization in Banach spaces
with uniformly convex penalty terms}, Numer. Math., 127 (2014), 485--513.

\bibitem{KN2008} S. Kindermann and A. Neubauer, {\it On the convergence of the quasioptimality criterion for
(iterated) Tikhonov regularization}, Inverse Probl. Imaging 2 (2008), no. 2, 291--299.

\bibitem{R1999} A. Rieder, {\it On the regularization of nonlinear ill-posed problems via inexact Newton
iterations}, Inverse Problems, 15 (1999), 309--327.

\bibitem{SEK1993} O. Scherzer, H. W. Engl and K. Kunisch, {\it Optimal a posteriori parameter choice
for Tikhonov regularization for solving nonlinear ill-posed problems}, SIAM J. Numer. Anal. 30 (1993), 1796--838.

\bibitem{SKHK2012} T. Schuster, B. Kaltenbacher, B. Hofmann and K. S. Kazimierski, {\it Regularization Methods
in Banach Spaces},  Radon Series on Computational and Applied Mathematics 10， Walter de Gruyter, Berlin 2012.

\bibitem{TJ2003} U. Tautenhahn and Q. Jin, {\it Tikhonov regularization and a posteriori rules
for solving nonlinear ill posed problems}, Inverse Problems, 19 (2003), no. 1, 1--21.

\bibitem{TGK1979} A. N. Tikhonov, V. B. Glasko and J. A. Kriksin, {\it On the question of quasi-optimal choice
of a regularized approximation}, Dokl. Akad. Nauk., 248 (1979), 531--535.

\bibitem{W1977} G. Wahba, {\it The approximate solution of linear operator equations when data are noisy},
SIAM J. Numer. Anal. 14 (1977), 651--667.


\bibitem{Z2002} C. Z\u{a}linscu, {\it Convex Analysis in General Vector Spaces}, World Scientific Publishing Co.,
Inc., River Edge, 2002.

\end{thebibliography}
\end{document}